\documentclass[12pt]{article}
\def\R{\bf R}  
\input epsf

\begin{document}

\centerline{\uppercase{\bf A short proof of the Kuratowski graph planarity criterion}
\footnote{This note is based on the author's lectures at the Kirov Region Summer School, St Petersburg Summer School, the Moscow Olympic School, mathematical circles at Kolmogorov College and at Moscow Center for Continuous Mathematical Education. I would like to acknowledge B. Mohar, D. Permyakov, V. Volkov, M. Vyalyi and T. Shaihieva for useful discussions.}
}
\bigskip
\centerline{\bf A. Skopenkov
\footnote{http://www.mccme.ru/\~\ skopenkov. Supported by Simons-IUM Foundation.}
}

\bigskip
\small
{\bf Abstract.} This paper is purely expositional.
The statement of the Kuratowski graph planarity criterion is simple and
well-known.
However, its classical proof is not easy.
In this paper we present the Makarychev proof (with further simplifications by
Prasolov, Telishev, Zaslavski and the author) which is possibly the simplest.
The paper is accessible for students familiar with the notion of a graph, and
could be an interesting easy reading for mature mathematicians.


\normalsize
\bigskip
A graph is called {\it planar} if it can be drawn in the plane without self-intersections.

\begin{figure}[h]
\centerline{\epsfbox{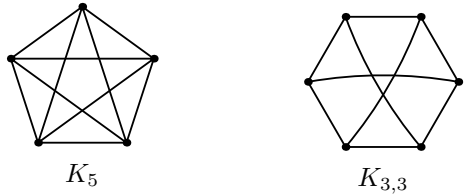}}
\caption{The Kuratowski graphs}
\label{k5k33}
 \end{figure}

\smallskip
{\bf The Kuratowski Theorem.}
{\it A graph is planar if and only if it does not contain a subgraph
homeomorphic to $K_5$ or to $K_{3,3}$ (fig. \ref{k5k33}).}

For definition of homeomorphic graphs, as well as for a short proof of the
`only if' part of the Kuratowski theorem see [Pr04].
For results related to this theorem see either the Russian version of this text or
[Cl34, Cl37, Ep66, GHV79, HJ64, Ku00, MS67, MA41, RS90, RS99, Sa91, Sk95, Sk08, Sk, Wh33].

Here we present a simple proof of the `if' part of the Kuratowski Theorem based on [Ma97], cf. [Th81], with further simplifications by Prasolov, Telishev, Zaslavski and the author.


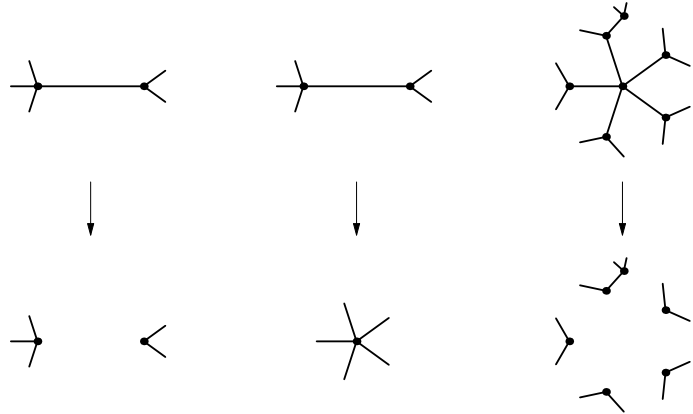
\begin{figure}[h]
\centerline{\epsfbox{graphs.4}}
\centerline{Deletion of an edge: $G-e$, contraction of an edge: $G/e$, and deletion of a vertex: $G-x$}
\label{udal}
\end{figure}

Clearly, it suffices to prove the Kuratowski Theorem for graphs without loops and multiple edges.
So we consider only such graphs.
By contraction of an edge we would understand contraction of an edge with replacement
of each obtained edge of multiplicity greater than 1 by an edge of multiplicity 1.

We prove the `if' part of the Kuratowski Theorem in the following equivalent form.

\smallskip
{\bf Proposition.}
{\it If a connected graph $G$ is not isomorphic to $K_5$ or to $K_{3,3}$, and for each edge $e$ of $G$ both graphs $G-e$ and $G/e$ are planar,
then $G$ is planar.}

\begin{figure}[h]
\centerline{\epsfbox{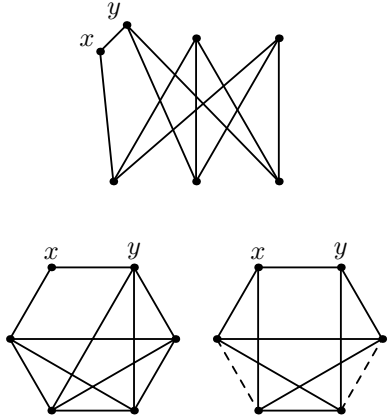}}
\caption{'Uncontraction of an edge' in the Kuratowski graphs}\label{rast}
\end{figure}

\smallskip
{\it Proof that Proposition implies the `if' part of the Kuratowski Theorem.}
 The statement `{\it graph $G$ contains a subgraph homeomorphic to the graph $H$}'
will be abbreviated to '$G\supset H$'.

The `if part' of the Kuratowski Theorem is proved by induction on the number of edges in the graph.
By Proposition the inductive step follows because
{\it $G\supset K_5$ or $G\supset K_{3,3}$ if either $G-e\supset K_5$ or $G-e\supset K_{3,3}$
or $G/e\supset K_5$ or $G/e\supset K_{3,3}$ for some edge $e$ of $G$}.
`For $G-e$' the italicized assertion is obvious.
If $G/e\supset K_{3,3}$ then $G\supset K_{3,3}$, and if $G/e\supset K_5$ then $G\supset K_5$ or $G\supset K_{3,3}$ (fig. \ref{rast}).
QED

\smallskip
{\bf Lemma on the Kuratowski Graphs.}
{\it For each graph $G$ the following three conditions are equivalent:

(1) For each edge $xy$ of $G$ the graph $G-x-y$ does not contain
$\theta$-graphs, and from each vertex of $G-x-y$ at least two edges are issuing.

(2) For each edge $xy$ of the graph $G$ the graph $G-x-y$ is a cycle
(containing $n\ge3$ vertices).

(3) $G$ is isomorphic either to $K_5$ or to $K_{3,3}$.}

\smallskip
The implications $(3)\Rightarrow(2)\Rightarrow(1)$ in the Lemma on the Kuratowski Graphs
are clear and are not used in the proof of the Kuratowski theorem.

\begin{figure}[h]
\centerline{\epsfbox{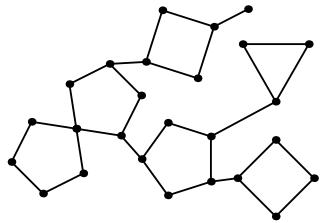}}
\caption{A 'tree' of cycles}\label{derevo}
\end{figure}

\smallskip
{\it Proof of the implication $(1)\Rightarrow(2)$ in the Lemma on the Kuratowski Graphs.}
Condition (1) implies that $K-x-y$ is a disjoint union of `trees' whose `vertices' are cycles.
(fig. \ref{derevo}; formally, each block of $K$ is a cycle).
Therefore $K-x-y$ contains a `hanging' cycle, i.e. a cycle $C$ having only one common
vertex $v$ with the remaining graph.
This cycle $C$ has at most two other vertices $p$ and $q$.
Since $K-x-y$ does not contain isolated vertices, from each vertex of $K$ there
issue at least three edges. Hence each of these vertices
$p$ and $q$ is joined either with $x$ or with $y$.
Therefore in the union of $C$ and the edges of $K$ joining vertices $x,y,p,q$
we can find a $\theta$-subgraph.
Hence by (1) each edge of $K-x-y$ has an end on $C$.
Since by (1) the graph $K-x-y$ does not contain hanging vertices, this graph coincides
with $C$.
QED

\smallskip
{\it Proof of the implication $(2)\Rightarrow(3)$ in the Lemma on the Kuratowski Graphs.}
If $n=3$ then for each two vertices $b$ and $c$ of the cycle $K-x-y$ the graph $K-b-c$
is a cycle, hence the remaining vertex of the cycle $K-x-y$ is joined (by an edge of $K$)
both to $x$ and to $y$.
Hence $K=K_5$.

If $n\ge4$, then take any four consecutive vertices $a,b,c,d$ of the cycle $K-x-y$.
Since $K-b-c$ is a cycle, in the graph $K$ one of the vertices $a$ or $d$ is
joined (by an edge) to $x$ (and not joined to $y$), the other is joined to $y$
(and not joined to $x$), whereas the vertices of the cycle $K-x-y$ different from $a,b,c,d$
(which do not exist when $n=4$) are not joined neither to $x$ nor to $y$.
For $n\ge5$ we obtain a contradiction.
For $n=4$ we see that the four vertices of the cycle $K-x-y$ are joined to $x$ and to $y$ one after the other,
hence $K=K_{3,3}$.
QED

\smallskip
{\it Proof of Proposition.} By the implication $(1)\Rightarrow(3)$ of Lemma on the Kuratowski Graphs
there is an edge $xy$ of $G$ such that $G-x-y$ contains either a vertex of degree at most 2 (in $G-x-y$) or $\theta$-graph.

If the degree of some vertex of $G$ is 1 or 2, then contraction of one of them gives a planar graph, so $G$ is planar.
So assume that in $G$ out of each vertex there issues at least three edges.

Hence graph $G-x-y$ does not have isolated vertices and if it has a hanging vertex $p$, then $p$ is joined both to $x$ and to $y$ in $G$.
Draw the graph $G-(xy)$ in the plane without self-intersections.
Add edge $xy$ along edges $px$ and $py$.
We obtain a drawing of $G$ in the plane without self-intersections.

\begin{figure}[h]
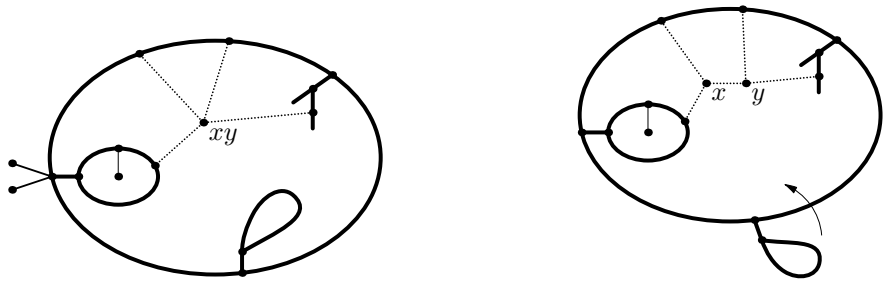

  \centerline{
    \begin{tabular}{c@{\hskip1in}c}
  \epsffile{graphs.6}&\epsffile{graphs.7}
  \end{tabular}}
\caption{Drawings of the graphs $G/xy$ and $G$ in the plane}\label{gg}
\end{figure}

Now consider the case when $G-x-y$ has an $\theta$-subgraph.
Draw graph $G/xy$ in the plane without self-intersections (fig. \ref{gg} left).
Drawing of graph $G-x-y=G/xy-xy$ in the plane is obtained by deleting
the edges of $G/xy$ issuing out of the vertex $xy$.
Take the face of (the image of) $G/xy-xy$ that contains the vertex $xy$ of the graph $G/xy$.
Denote by $\overline C$ the boundary of this face.

Observe that {\it the boundary of a face cannot contain a $\theta$-subgraph.}

(This statement could be derived from the Jordan Curve Theorem.
Another proof could be obtained supposing the contrary.
If the boundary of a face contains a $\theta$-subgraph, then take a point inside
this face and join it by three edges to three points on three `arcs'
of the $\theta$-subgraph.
We obtain a drawing of $K_{3,3}$ in the plane without self-intersections.
Contradiction.)

So $G-x-y\ne\overline C$.
Then edges of $G-x-y-\overline C$
\footnote{{\it Deletion of a subgraph} is deletion of all the edges of this subgraph and of all the vertices which are endpoints
only of edges of the subgraph.
Note that deletion of a vertex is not the same as deletion of a subgraph formed by this vertex. }
are contained in a face of (the image of) $G/xy-xy$ not containing the vertex $xy$.
Hence graph $\overline C$ splits the plane.
Therefore there is a cycle $C\subset \overline C$ such that $xy$ is (without loss of generality)
inside $C$ and certain edge of $G-x-y-\overline C$ is outside $C$.

Denote by $R$ the union of all the edges of $G/xy$ lying outside $C$.
(Possibly $G-x-y-\overline C\ne R$.)
We may assume that $R$ is a subgraph of $G$.
Draw graph $G-R$ in the plane without self-intersections (fig. \ref{gg} right).
We may assume that in this drawing the edges of $G$ issuing out of $x$ or $y$ are inside the cycle $C$.

Each connected component of $G-x-y-R-C$ intersects $C$ at most by one point.

(Indeed, otherwise $G-x-y-R-C$ has a path joining two points of $C$.
In the drawing of $G/xy$ the corresponding path lies inside $C$.
Hence the path splits the interior of $C$ into two parts, one of the containing $xy$, and the other is not contained in the face bounded
by $\overline C$. Hence $C\not\subset\overline C$, which is a contradiction.)

Therefore we can shift to the interior of  $C$ each connected component of $G-x-y-R-C$ (see an arrow in fig. \ref{gg}).
 Then $G-R-C$ is inside $C$.
Draw $R$ outside $C$ like for the drawing of $G/xy$ (fig. \ref{gg} left).
We obtain a drawing of $G$ in the plane without self-intersections.
QED

\small
\bigskip
\centerline{\bf References}
\smallskip


[Cl34] S.~Claytor, Topological immersions of peanian continua in a spherical
surface, Ann. of Math. 35 (1934), 809--835.

[Cl37] S.~Claytor, Peanian continua not embeddable in a spherical surface,
Ann. of Math. 38 (1937) 631--646.

[Ep66] D. B. A. Epstein, Curves on 2-manifolds and isotopies, Acta Math. 115
(1966) 83--107.

[GHW79] H.~H.~Glover, J.~P.~Huneke and C.~S.~Wang, 103 graphs that are
irreducible for the projective plane, J.~Comb.~Th., 27:3 (1979) 332--370.

[HJ64] R.~Halin and H.~A.~Jung, Karakterisierung der komplexe der Ebene und der
2-Sph\"are, Arch. Math. 15 (1964) 466--469.

[Ku00] V. A. Kurlin, Basic embeddings into products of graphs, Topol. Appl.
102 (2000) 113--137.

[Ma97] Yu. Makarychev, A short proof of Kuratowski's graph planarity criterion,
J. of Graph Theory, 25 (1997) 129--131.

[MS67] S.  Marde\v si\' c and J.  Segal, $\varepsilon$-mappings and generalized
manifolds, Michigan Math. J. 14 (1967) 171--182.

[MA41] S.~McLane and V.~W.~Adkisson, Extensions of homeomorphisms on the
spheres, Michig. Lect. Topol., Ann~Arbor, (1941) 223--230.

[Pr04] V. Prasolov, Elements of combinatoril and differential topology, AMS translations.

[RS99] D. Repovs and A. Skopenkov, New results on embeddings of polyhedra
and manifolds into Euclidean spaces (in Russian), Uspekhi Mat. Nauk, 54:6
(1999) 61--109. English transl.: Russ. Math. Surv. 54:6 (1999),
1149--1196.

[RS90] N.~Robertson and P.~D.~Seymour, Graph minors VIII, A Kuratowski graph
theorem for general surfaces, J.~Comb. Theory, 48B (1990) 255--288.

[Sa91] K.~S.~Sarkaria, Kuratowski complexes, Topology, 30 (1991) 67--76.

[Sk95] A.Skopenkov, A description of continua basically embeddable in $\R^2$,
Topol. Appl. 65 (1995) 29--48.


[Sk08] A. Skopenkov, Embedding and knotting of manifolds in Euclidean spaces, in:
Surveys in Contemporary Mathematics, Ed. N. Young and Y. Choi
London Math. Soc. Lect. Notes, 347 (2008) 248--342. arxiv:math/0604045

[Sk] A. Skopenkov, Algebraic topology from an elementary viewpoint, in Russian,
MCCME, Moscow, to appear. arXiv:0808.1395

[Th81] C.~Thomassen, Kuratowski's theorem, J.~Graph. Theory, 5 (1981) 225--242.

[Wh33] H. Whitney, Planar graphs, Fund. Math. 21 (1933) 73--84.

\end{document}


\centerline{\uppercase{\bf Вокруг критерия Куратовского планарности графов}
\footnote{Эта заметка является дополненной версией заметки [Sk05]. 
Она основана на занятиях, проведенных в А. Скопенковым в Кировской ЛМШ, Московской ОВШ и на кружках `Математический семинар', `Олимпиады и математика'.
Благодарю В. Волкова, М. Вялого, Д. Пермякова и Т. Шайхиеву  за полезные замечания и обсуждения, М. Вялого и издательство
МЦНМО за подготовку рисунков, а также Б. Мохара и С. Матвеева за предоставленные ссылки.}
}
\bigskip
\centerline{\bf А.Б.~Скопенков
\footnote{http://www.mccme.ru/\~\ skopenkov. Поддержан грантом фонда Саймонса.}
 }

\bigskip
\small
{\bf Аннотация.}
Формулировка критерия Куратовского планарности графов хорошо известна
(все необходимые понятия и эта формулировка напомнены в тексте).
Однако его классическое доказательство сложно и приводится не во всех книгах по теории графов.
Здесь приводится доказательство Макарычева (с упрощениями, сделанными Заславским, Прасоловым, Телишевым и автором);
по-видимому, оно является наиболее простым.
Перед доказательством критерия Куратовского приводятся необходимые определения,
а после --- близкие результаты.
\normalsize

\bigskip
\centerline{\bf Планарные графы}
\smallskip

Граф называется {\it планарным}, если его можно без самопересечений нарисовать на плоскости.
Например, любое дерево и любой граф, образованный вершинами и ребрами
некоторого выпуклого многогранника --- планарные.
А графы $K_5$ и $K_{3,3}$ (см. рис.) не являются планарными.
Это можно доказать путем небольшого перебора с использованием следующей теоремы [Pr04, \S1, Теорема 1.3].
(Доказательство непланарности графа $K_5$, основанное на понятии коэффициента пересечения, см. в [BE82, Sk05, Pr04, \S1].)

\begin{figure}[h]
\centerline{\epsfbox{graphs.2}}
\caption{Графы Куратовского}\label{k5k33}
 \end{figure}

\smallskip
{\bf Теорема Жордана.}
{\it Замкнутая несамопересекающаяся кривая (т.е. цикл) делит плоскость ровно на две части.
(При этом одна часть ограничена, другая неограничена, причем две точки
плоскости, не принадлежащие кривой, лежат в одной части тогда и только тогда,
когда их можно соединить ломаной, не пересекающей кривой).}

\smallskip
Обсуждение и доказательство этой теоремы см., например, в [An03].

{\it Плоским графом} называется изображение графа на плоскости без
самопересечений.
Иногда такое изображение называют просто графом, но это неточно, поскольку
планарный граф можно изобразить (без самопересечений) на плоскости разными способами.

\smallskip
\centerline{\epsfbox{pict1.2}}
\centerline{Различные изображения графа на плоскости}


\smallskip
Грубо говоря, {\it подграф} данного графа --- это его часть.
Формально, граф $G$ называется {\it подграфом} графа $H$, если множество
вершин графа $G$ содержится в множестве вершин графа $H$ и каждое ребро
графа $G$ является ребром графа $H$.
При этом две вершины графа $G$, соединенные ребром в графе $H$, не
обязательно соединены ребром в графе $G$.

Ясно, что любой подграф планарного графа планарен.

 Грубо говоря, графы изоморфны, если они одинаковы (при этом их изображения
на плоскости могут быть разными).
Формально, графы $G_1$ и $G_2$ (без петель и краьных ребер) называются
{\it изоморфными}, если существует
взаимно-однозначное отображение $f:V_1\to V_2$ множества $V_1$ вершин графа
$G_1$ на множество $V_2$ вершин графа $G_2$, удовлетворяющее условию:
{\it вершины $A,B\in V_1$ соединены ребром в том и только в том случае, если
вершины $f(A),f(B)\in V_2$ соединены ребром}.


Операция {\it подразделения ребра} графа показана на рисунке.

\smallskip
\centerline{\epsfbox{graphs.1}}
\centerline{Подразделение ребра}

\smallskip
Два графа называются {\it гомеоморфными}, если от одного можно перейти к
другому при помощи операций подразделения ребра и обратных к ним.
Или, эквивалентно, если существует граф $G$, полученный из обоих данных графов
операциями подразделения ребра.

Ясно, что гомеоморфные графы являются или не являются планарными одновременно.

\smallskip
{\bf Теорема Куратовского.} {\it Граф является планарным тогда и только тогда,
когда он не содержит подграфа, гомеоморфного графу $K_5$ или $K_{3,3}$ (рис. k5k33).}

\smallskip
Эта теорема объявлена также замечательным советским математиком Львом
Семеновичем Понтрягиным (доказательство не опубликовано), а также Фринком
и Смитом.
Поэтому иногда ее называют теоремой Понтрягина-Куратовского.
В 1920-е годы Карл Менгер объявил, что {\it граф, степень каждой вершины
которого равна 3, является плоским тогда и только тогда, когда он не
содержит подграфа, гомеоморфного графу $K_{3,3}$}.
Читатель может попытаться самостоятельно доказать этот факт (вытекающий из
теоремы Куратовского).
Кроме теоремы Куратовского, существует много других критериев планарности
графов [Th81].
Огромный интерес к поиску критерия планарности графов объясняется, в частности,
наличием одной из величайших математических гипотез --- гипотезы четырех красок
[Pr04, \S1].
Она утверждает, что вершины любого плоского графа можно правильно раскрасить
в 4~цвета.
Раскраска вершин (граней) плоского графа называется {\it правильной},
если любые две соседние вершины (грани) окрашены в разные цвета.


\bigskip
\centerline{\bf Простое доказательство достаточности в теореме Куратовского}

Приводимое доказательство в основном принадлежите Ю. Макарычеву (он придумал свое доказательство, еще будучи школьником!), ср. [Th81, \S5].

\begin{figure}[h]
\centerline{\epsfbox{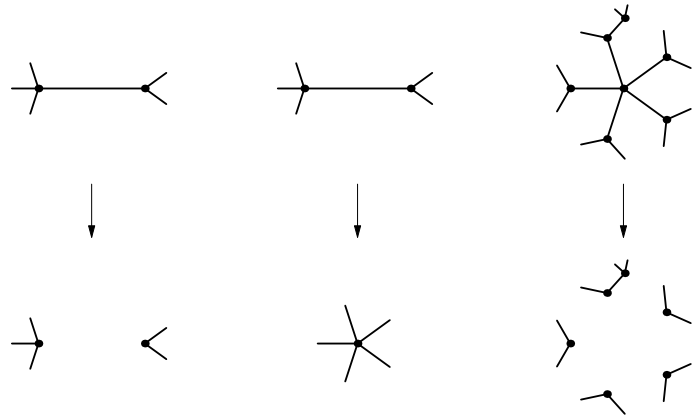}}
\caption{Удаление ребра $G-e$, стягивание ребра $G/e$ и удаление вершины $G-x$}
\label{udal}
\end{figure}

Так как достаточность в теореме Куратовского достаточно доказать для графов без петель и кратных ребер,
то будем рассматривать только такие графы.
Под стягиванием ребра будем понимать стягивание ребра вместе с заменой
каждого получившегося ребра кратности больше 1 на ребро кратности 1.

\smallskip
{\bf Утверждение.}
{\it Если связный граф $G$
не изоморфен ни $K_5$, ни $K_{3,3}$, и для любого ребра $e$ графа $G$ оба графа $G-e$ и $G/e$ планарны, то $G$ планарен.}

\begin{figure}[h]
\centerline{\epsfbox{graphs.5}}
\caption{'Растягивание ребра' в графах Куратовского}\label{rast}
\end{figure}

\smallskip
{\it Доказательство достаточности в теореме Куратовского с использованием Утверждения.}
Свойство `{\it граф $G$ содержит подграф, гомеоморфный графу $H$}'
будем сокращенно записывать в виде `$G\supset H$'.

Достаточность в теореме Куратовского доказывается индукцией по количеству ребер в графе.
Шаг индукции следует из Утверждения, поскольку
{\it  если $G-e\supset K_5$ или $G-e\supset K_{3,3}$ или $G/e\supset K_5$ или $G/e\supset K_{3,3}$ для некоторого ребра $e$ графа $G$,
то $G\supset K_5$ или $G\supset K_{3,3}$}.
`Для $G-e$' выделенное курсивом утверждение очевидно.
Если $G/xy\supset K_{3,3}$, то $G\supset K_{3,3}$, а если
$G/xy\supset K_5$, то $G\supset K_5$ или $G\supset K_{3,3}$ (рис. \ref{rast}).
 QED

\smallskip
{\bf Лемма о графах Куратовского.}
{\it Для произвольного графа $K$ следующие три условия
равносильны:

(1) Для любого ребра $xy$ графа $K$ граф $K-x-y$ не содержит
$\theta$-графа, и из каждой вершины графа $K-x-y$ выходит не менее двух ребер.

(2) Для любого ребра $xy$ графа $K$ граф $K-x-y$ является циклом
(содержащим $n\ge3$ вершин).

(3) $K$ изоморфен $K_5$ или $K_{3,3}$.}

\smallskip
 Импликации $(3)\Rightarrow(2)\Rightarrow(1)$ в лемме о графах Куратовского
очевидны и не используются в доказательстве теоремы Куратовского.

\begin{figure}[h]
\centerline{\epsfbox{graphs.9}}
\caption{'Дерево' из циклов}\label{derevo}
\end{figure}

\smallskip
{\it Доказательство импликации $(1)\Rightarrow(2)$ в лемме о графах Куратовского.}
Ввиду (1) в графе $K-x-y$ существует `висячий' цикл, т.е. цикл $C$, имеющий с остальным графом только одну общую вершину $v$
(ибо граф $K-x-y$ представляет собой одно или несколько `деревьев', `вершинами' которых служат циклы, рис. \ref{derevo}; формально говоря, каждым блоком графа $K-x-y$ является цикл).
В этом цикле $C$ есть еще по крайней мере две вершины $p$ и $q$.
Так как в графе $K$ нет вершин, из которых выходит менее трех ребер, то каждая
из этих вершин $p$ и $q$ соединена либо с $x$, либо с $y$.
Поэтому в объединении цикла $C$ и ребер графа $K$, соединяющих вершины $x,y,p,q$,
можно выделить $\theta$-подграф.
Значит, по (1) каждое ребро графа $K-x-y$ имеет конец на цикле $C$.
Поскольку по (1) граф $K-x-y$ не содержит висячих вершин, то $K-x-y=C$.
QED

\smallskip
{\it Доказательство импликации $(2)\Rightarrow(3)$ в лемме о графах Куратовского.}
При $n=3$ для любых двух вершин $b$ и $c$ цикла $K-x-y$ граф $K-b-c$
является циклом, поэтому оставшаяся вершина цикла $K-x-y$ соединена (ребром)
в $K$ и с $x$, и с $y$.
Поэтому $K=K_5$.

При $n\ge4$ возьмем любые четыре соседние вершины $a,b,c,d$ цикла $K-x-y$.
Поскольку граф $K-b-c$ является циклом, то в $K$ одна из вершин $a$ и $d$
соединена с $x$ (и не соединена с $y$), другая соединена с $y$ (и не соединена
с $x$), а отличные от $a,b,c,d$ вершины цикла $K-x-y$ (которых нет при $n=4$)
не соединены ни с $x$, ни с $y$.
При $n\ge5$ получаем противоречие.
При $n=4$ получаем, что четыре вершины цикла $K-x-y$ соединены с $x$ и $y$
попеременно, откуда $K=K_{3,3}$.
QED

\smallskip
{\it Доказательство утверждения.}
 Так как $G$ не изоморфен ни $K_5$, ни $K_{3,3}$, то по лемме о графах Куратовского существует ребро $e=(xy)$ графа $G$, для которого
в графе $G-x-y$ найдется либо  вершина степени меньше 2 (в $G-x-y$), либо $\theta$-подграф.

Если в графе $G$ из некоторой вершины выходит одно или два его ребра, то при стягивании одного из них получается планарный граф,
значит, и граф $G$ планарен.
Поэтому далее будем считать, что из каждой вершины графа $G$ выходит не менее трех его ребер.

Поэтому в графе $G-x-y$ нет изолированных вершин, и если есть висячая вершина $p$, то она соединена и с $x$, и с $y$ в графе $G$.
Нарисуем граф $G-(xy)$ на плоскости без самопересечений.
Так как в графе $G$ из $p$ выходит три ребра, то `с одной стороны' от пути $xpy$ из $p$ не выходит ребер.
'Подрисуем' ребро $xy$ вдоль пути $xpy$ `с этой стороны' от него.
Получим изображение графа $G$ на плоскости без самопересечений.

Рассмотрим теперь случай, когда в графе $G-x-y$ найдется $\theta$-подграф.

Из теоремы Жордана следует, что {\it любой плоский граф разбивает плоскость на
конечное число связных частей.}
Эти части называются {\it гранями} плоского графа.

\begin{figure}[h]
  \centerline{
    \begin{tabular}{c@{\hskip1in}c}
  \epsffile{graphs.6}&\epsffile{graphs.7}
  \end{tabular}}
\caption{Изображение на плоскости графов $G/xy$ и $G$}\label{gg}
\end{figure}

Нарисуем без самопересечений на плоскости  граф $G/xy$ (рис. \ref{gg} слева).
Изображение графа $G-x-y=G/xy-xy$ на плоскости получается стиранием
ребер графа $G/xy$, выходящих из вершины $xy$.
Обозначим через $\overline C$ границу той грани (изображения) графа $G/xy-xy$,
которая содержит вершину $xy$ графа $G/xy$.

Заметим, что {\it граница грани не может содержать $\theta$-подграфа.}

(Это утверждение можно вывести из теоремы Жордана.
Другое доказательство получается от противного: если граница грани содержит
$\theta$-подграф, то возьмем точку внутри этой грани и соединим ее тремя
ребрами с тремя точками на трех 'дугах' $\theta$-подграфа.
Получим изображение графа $K_{3,3}$ на плоскости без самопересечений.
Противоречие.)

Поэтому $G-x-y\ne\overline C$.
Тогда ребра графа $G-x-y-\overline C$
\footnote{Удаление подграфа --- удаление всех его ребер и всех вершин, из которых выходят только ребра этого подграфа.
Заметим, что удаление вершины --- не то же самое, что удаление подграфа из этой вершины.}
находятся в грани (изображения) графа $G-x-y=G/xy-xy$, не содержащей вершины $xy$.
Значит, граф $\overline C$ разбивает плоскость.
Поэтому найдется цикл $C\subset\overline C$, относительно которого вершина $xy$
лежит (не уменьшая общности) внутри, а некоторое ребро графа $G-x-y-\overline C$ --- вне.

Обозначим через $R$ объединение всех ребер графа $G/xy$, лежащих вне цикла $C$.
(Возможно, $G-x-y-\overline C\ne R$.)
Можно считать, что $R$ --- подграф в $G$ (а не только в $G-x-y$).

Граф $G-R$
можно нарисовать на плоскости без самопересечений (рис. \ref{gg} справа).
Можно считать, что ребра графа $G$, выходящие из $x$ или $y$, на изображении графа $G-R$ лежат {\it внутри} цикла $C$.

Каждая компонента связности графа $G-x-y-R-C$ пересекается с $C$ не более чем по одной точке.

(Если это не так, то в $G-x-y-R-C$ есть путь, соединяющий две точки на $C$.
На изображении графа $G/xy$ соответствующий путь лежит внутри цикла $C$.
Значит, этот путь разбивает внутреннюю часть цикла $C$ на две части, одна из которых содержит $xy$, a другая не лежит в грани,
ограниченной $\overline C$. Поэтому $C\not\subset\overline C$ --- противоречие.)

Поэтому можно перекинуть внутрь цикла $C$ каждую компоненту связности графа $G-x-y-R-C$ (см. стрелочку на рис. \ref{gg} справа).
Значит, граф $G-R-C$ можно нарисовать внутри цикла $C$.
 Нарисуем $R$ вне $C$, как для изображения графа $G/xy$ (рис. \ref{gg} слева).
Получим изображение графа $G$ на плоскости без самопересечений.
QED



\bigskip
\centerline{\bf Запрещенные подсистемы}

\smallskip
Если некоторая подсистема системы~$N$ не реализуема в~другой системе $M$, то и
$N$ не реализуема в~$M$.
Естественная идея --- попытаться найти список 'запрещенных' систем
$N_1,\dots,N_k$, не реализуемых в $M$, со следующим свойством:

{\it для того, чтобы система $N$ была реализуема в~$M$ необходимо и достаточно,
чтобы $N$ не содержал ни одной из этих 'запрещенных' подсистем.}

Классический пример теоремы такого рода --- теорема Куратовского.
Так можно описать графы, вложимые в любую данную поверхность [RS90], а также
много других классов графов
или более общих обьектов (например, графы и даже пеановские континуумы,
{\it базисно} вложимые в плоскость [Sk95, Ku00]).
Заметим, что список запрещенных подграфов для вложимости графа в лист
Мебиуса содержит целых 103 графа [GHW79].
Даже {\it существование} такого конечного списка для произвольной поверхности
доказывется сложно [AH89, RS90].
Список запрещенных полиэдров бесконечен для вложимости двумерных
полиэдров в~$\R^3$ или $n$-мерных полиэдров в~$\R^{2n}$, где $n\ge2$ [Sa91].
Поэтому интересны другие препятствия к вложимости.
Заметим, что одно из самых полезных препятствий строится с помощью
{\it конфигурационного пространства} упорядоченных пар различных точек данного
пространства [RS99, Sk08].

Приведем формулировки некоторых результатов (доказательства оставляем читателю в качестве задач).
Те формулировки, в которых встречаются неизвестные читателю объекты, он может
игнорировать.

\smallskip
{\bf Теорема Шартрана-Харари.} {\it Граф $G$ можно нарисовать на плоскости без
самопересечений так, чтобы он был границей некоторой} одной {\it грани тогда и
только тогда, когда $G$ не содержит $\theta$-подграфа.}

\smallskip
Назовем несамопересекающийся цикл $C$ в связном графе $G$ {\it граничным},
если существует изображение без самопересечений графа $G$ на плоскости, при
котором цикл $C$ изображается границей некоторой грани.
Следующий результат можно вывести из теоремы Куратовского, ср. [Cl34].

\begin{figure}[h]
  \centerline{
    \begin{tabular}{c@{\hskip1in}c}
  \epsffile{graphs.13}&\epsffile{graphs.14}
  \end{tabular}
  }
\caption{Цикл $C$ не может быть границей внешней грани}\label{outer}
\end{figure}

{\bf Относительная версия теоремы Куратовского.}
{\it Цикл $C$ является граничным тогда и только тогда, когда граф $G$ планарен
и цикл $C$ не содержится в подграфе графа $G$, как на рис.~\ref{outer}.}

\smallskip
А вот следующий результат проще доказывать, не используя теорему Куратовского
(подробнее см. [Sk, глава 1]).

\smallskip
{\bf Теорема о 8 и $\theta$.}
{\it Граф $G$ с заданными циклами ребер, выходящих из каждой вершины, можно так
изобразить без самопересечений на плоскости, чтобы указанные циклы шли по
часовой стрелке, тогда и только тогда, когда $G$ не содержит 'восьмерки' или
'буквы $\theta$' с циклами, изображенными на рис.~\ref{cycle}.}

\begin{figure}[h]
  \centerline{\epsffile{graphs.15}  }
\caption{Графы с циклическими порядками, не реализуемыми на плоскости}
\label{cycle}
\end{figure}


Два вложения (т.е. изображения без самопересечений)  $f,g$ одного и
того же графа в плоскость называются {\it изотопными}, если одно можно
так непрерывно продеформировать в другое, чтобы в процессе деформации мы все
время имели бы вложение (формальное определение см., например, в [Ma07]).

\smallskip
{\bf Теорема Маклейна-Эдкиссона.}
{\it Два вложения связного графа в плоскость изотопны тогда и только тогда,
когда их сужения на любой триод $T$ и на любой несамопересекающийся цикл $S^1$
изотопны (т.е. не таковы, как на рис.~\ref{mcle})} [MA41].

\begin{figure}[h]
  \centerline{
    \begin{tabular}{c@{\hskip1in}c}
  \epsffile{graphs.16}&\epsffile{graphs.17}
  \end{tabular}}
\caption{Различные вложения триода и окружности в плоскость}\label{mcle}
\end{figure}

Эту теорему удобно сначала доказать для деревьев, а потом свести общий случай
к случаю деревьев путем выделения максимального дерева.

Теорема Маклейна-Эдкиссона справедлива также для {\it полиэдра} или даже
{\it пеановского континуума} [MA41].

Теорема Маклейна-Эдкиссона (без утверждения в скобках) справедлива для вложений
в сферу, тор и другие сферы с ручками (доказательство аналогично).
Заметим, что любая изотопия графа на поверхности объемлема [С. В. Матвеев,
частное сообщение].

\smallskip
{\bf Теорема Баэра-Эпштейна.}
{\it Две замкнутые несамопересекающиеся кривые на двумерном многообразии
гомотопны
тогда и только тогда, когда они изотопны} [Ep66].

\smallskip
Теорема Эпштейна сводит вопрос о классификации вложений окружности в
двумерное многообразие $N$ (и, тем самым, произвольного графа в сферу с
ручками и дырками) к вопросу о реализуемости элементов фундаментальной группы
$\pi_1(N)$ вложенными окружностями. Но последний вопрос очень сложен.

\smallskip
{\bf Задача.}
Граф называется (вершинно) {\it $k$-связным}, если он остается связным после
удаления любой $k-1$ вершины и распадается после удаления некоторых $k$ вершин.
Выведите из теоремы Маклейна-Эдкиссона следующие утверждения.

Любое вложение произвольного трехсвязного графа в сферу может быть получено из
любого другого композицией изотопии и осевых симметрий.

Любое вложение двусвязного графа в сферу может быть получено из любого
другого композицией изотопии и 'переворачиваний блоков' (рис.~\ref{block}).

\begin{figure}[h]
  \centerline{\epsffile{graphs.18}  }
\caption{Переворачивание блока}
\label{block}
\end{figure}

Определите операции, при помощи которых можно получить любое вложение
1-связного ($\Leftrightarrow$ связного) графа в сферу из любого другого.
Сделайте то же и для 0-связного ($\Leftrightarrow$ произвольного) графа.
Таким образом получится другое описание вложений
графов в плоскость с точностью до изотопии [Wh33].

\bigskip
\centerline{\bf Приложение: планарность полиэдров и континуумов}

\smallskip
Полиэдр (синоним: тело симплициального комплекса) --- это многомерный аналог
графа.
Определение см., например, в [Pr04, \S8, Ma07].
Уже двумерные полиэдры --- интересные и сложные объекты, про которые имеется
несколько знаменитых и трудных нерешенных проблем [Ma07].
Поэтому удивительно, что имеется следующий результат.

\smallskip
{\bf Теорема Халина-Юнга.} {\it Связный полиэдр вложим в~$S^2$ тогда и
только тогда, когда он не содержит графов $K_5$, $K_{3,3}$ или
'зонтика'~$U^2$ (рис.~\ref{U})} [HJ64, MS67].

\begin{figure}[h]
  \centerline{\epsffile{graphs.19}  }
\caption{Зонтик}
\label{U}
\end{figure}

В этом результате интересна лишь часть 'тогда' и лишь для двумерных полиэдров.
Следующее доказательство (видимо, являющееся фольклорным) проще представленного
в [HJ64] и тем более в [MS67].

\smallskip
{\it Набросок доказательства части 'тогда'.}
Пусть связный 2-полиэдр $N\not\cong S^2$ не содержит ни графов $K_5$, $K_{3,3}$,
ни зонтика $U$.
Так как $N$ не содержит зонтика, то окрестность любой точки в $N$ является
объединением дисков и отрезков, склеенных за одну точку (рис. ~\ref{veer}
слева).
Если этих дисков больше одного, то заменим эту окрестность на изображенную на
рис.~\ref{veer} справа.

\begin{figure}[h]
  \centerline{
    \begin{tabular}{c@{\hskip1in}c}
  \epsffile{changes.1}&\epsffile{changes.2}
  \end{tabular}
  }
\caption{Преобразование окрестности точки}\label{veer}
\end{figure}

При этом преобразовании не появится подграфов $K_5$ и $K_{3,3}$;
Обратное преобразование является стягиванием 'звезды с несколькими лучами' и
поэтому сохраняет планарность.
Значит, достаточно доказать теорему для полученного 2-полиэдра.
Рассмотрим объединение $\bar N$ его двумерных граней.
Тогда окрестность любой точки в $\bar N$ является диском.
Значит, по теореме классификации поверхностей $\bar N$ является сферой с
ручками, пленками Мебиуса и дырками.
Поскольку каждый из графов $K_5$ и $K_{3,3}$ вложим и в тор с дыркой,
и в лист Мебиуса, то $\bar N$ есть несвязное объединение дисков с дырками.
Заменим каждый из этих дисков с дырками на граф с рис.~\ref{koleso}.
Полученный граф планарен.
По вложению этого графа в плоскость легко построить вложение полиэдра $N$
в плоскость. QED

\begin{figure}[h]
  \centerline{\epsffile{changes.3}}
\caption{Преобразование диска с дырками}\label{koleso}
\end{figure}

\smallskip
В терминах запрещенных подсистем можно также описать 'компактно бесконечные
графы' (т.е. локально связные континуумы), вложимые в плоскость.
{\it Континуум} --- компактное связное метрическое пространство.
Континуумы
естественно появляются при изучении динамических систем (даже гладких!).
Континуум называется {\it локально связным} (или континуумом Пеано), если
для любой его точки $x$ и ее окрестности $U$ существует такая меньшая
окрестность $V$ точки $x$, что любые две точки из $V$ соединяются некоторым
путем, целиком лежащим в $U$
(или, эквивалентно, если он является непрерывным образом дуги).
Локально связные континуумы могут быть очень сложно устроены.
Поэтому удивительно, что имеется следующий результат.

\smallskip
{\bf Теорема Клэйтора.} {\it Пеановский континуум вложим в $S^2$
тогда и только тогда, когда он не содержит континуумов $K_5$,~$K_{3,3}$,
$C_{K_5}$ и $C_{K_{3,3}}$ (рис.~\ref{clay})} [Cl34, Cl37].

\begin{figure}[h]
  \centerline{\epsffile{graphs.22}}
\caption{Континуумы Куратовского}
\label{clay}
\end{figure}

\smallskip
{\it Построение континуумов $C_{K_5}$ и $C_{K_{3,3}}$.}
Возьмем ребро $ab$ графа $K_5$ и отметим на нем новую вершину $a'$.
Пусть $P=K_5-(aa')$.
Пусть $P_n$ копия графа $P$.
Обозначим через $a_n$ и $a'_n$ вершины графа $P_n$, соответствующие $a$ и $a'$.
Положим
$$C_{K_5}:=(P_1\bigcup\limits_{a'_1=a_2}P_2\bigcup\limits_{a'_2=a_3}P_3\dots)
\bigcup\limits_{x=0}I,$$
где $\{P_n\}$ --- последовательность графов на плоскости со стремящимися к
нулю диаметрами, сходящаяся к точке $x\notin\sqcup^\infty_{n=1}P_n$.
Точно так же можно определить континуум $C_{K_{3,3}}$, взяв в начале
$K_{3,3}$ вместо $K_5$.

\smallskip
{\it Доказательство невложимости в теореме Клэйтора.}
Докажем невложимость $C_{K_5}$ (доказательство невложимости $C_{K_{3,3}}$
аналогично).
Пусть, напротив, $f:C_{K_5}\to\R^2$ --- вложение.
Обозначим через $S:=P-a-a'$ окружность в $P$, составленную из ребер, не
содержащих вершин $a$ и $a'$.
Аналогично определим $S_n\subset P_n$.

Так как $S_n$ сходится к $x=0$, то $f1$ лежит вне $fS_n$ для достаточно
большого $n$.

Так как $fI$ --- путь между $f0$ и $f1$, лежащий вне $fS_n$, то $f0$ лежит вне
$fS_n$.

Так как $S_n$ сходится к $x=0$, то $fS_m$ лежит вне $fS_n$ и
$fS_l$ лежит вне $fS_m$ для достаточно больших $m<l$.
Но тогда $fa_m$ и $fa'_m$ лежат вне $fS_m$.
Это противоречит тому, что для любого вложения $g:P\to\R^2$ точки $ga$ и $ga'$
лежат по разные стороны от образа $gS$.
QED

\smallskip
Отметим, что в теореме Куратовского можно заменить плоскость $\R^2$ на сферу
$S^2$.
В теоремах  Халина-Юнга и Клэйтора можно заменить $S^2$ на $\R^2$, только
добавив запрещенный подполиэдр или подконтинуум $S^2$.

\bigskip

\small
{\bf Литература}
\smallskip

[An03] Д. В. Аносов, Отображения окружности, векторные поля и их применения,
МЦНМО, Москва, 2003.

[AH89] D.~Archdeacon and P.~Huneke, A Kuratowski theorem for non-orientable
surfaces, J. Comb. Th., Ser.~B, 46 (1989) 173--231.

[BE82] В. Г. Болтянский и В. А. Ефремович, Наглядная топология,
Наука, Москва, 1982.

[Cl34] S.~Claytor, Topological immersions of peanian continua in a spherical
surface, Ann. of Math. 35 (1934), 809--835.

[Cl37] S.~Claytor, Peanian continua not embeddable in a spherical surface,
Ann. of Math. 38 (1937) 631--646.

[Ep66] D. B. A. Epstein, Curves on 2-manifolds and isotopies, Acta Math. 115
(1966) 83--107.

[GHW79] H.~H.~Glover, J.~P.~Huneke and C.~S.~Wang, 103 graphs that are
irreducible for the projective plane, J.~Comb.~Th., 27:3 (1979) 332--370.

[HJ64] R.~Halin and H.~A.~Jung, Karakterisierung der komplexe der Ebene und der
2-Sph\"are, Arch. Math. 15 (1964) 466--469.


[Ku00] V. A. Kurlin, Basic embeddings into products of graphs, Topol. Appl.
102 (2000) 113--137.

[Ma07] С. В. Матвеев, Алгоритмическая топология и классификация трехмерных
многообразий, МЦНМО, Москва, 2007.

[Ma97] Yu. Makarychev, A short proof of Kuratowski's graph planarity criterion,
J. of Graph Theory, 25 (1997) 129--131.

[MS67] S.  Marde\v si\' c and J.  Segal, $\varepsilon$-mappings and generalized
manifolds, Michigan Math. J. 14 (1967) 171--182.

[MA41] S.~McLane and V.~W.~Adkisson, Extensions of homeomorphisms on the
spheres, Michig. Lect. Topol., Ann~Arbor, (1941) 223--230.

[Pr04] В. В. Прасолов, Элементы комбинаторной и дифференциальной топологии,
МЦНМО, Москва, 2004.

[RS99] Д. Реповш и А.  Скопенков, Новые результаты о вложениях полиэдров и
многообразий в евклидовы пространств, УМН 54:6 (1999) 61--109.

[RS90] N.~Robertson and P.~D.~Seymour, Graph minors VIII, A Kuratowski graph
theorem for general surfaces, J.~Comb. Theory, 48B (1990) 255--288.

[Sa91] K.~S.~Sarkaria, Kuratowski complexes, Topology, 30 (1991) 67--76.

[Sk95] A.Skopenkov, A description of continua basically embeddable in $\R^2$,
Topol. Appl. 65 (1995) 29--48.

[Sk05] А. Скопенков, Вокруг критерия Куратовского планарности графов,
Мат. Просвещение, 9 (2005), 116--128 и 10 (2006), 276--277.
http://www.mccme.ru/free-books/matprosa.html

[Sk08] A. Skopenkov, Embedding and knotting of manifolds in Euclidean spaces, in:
Surveys in Contemporary Mathematics, Ed. N. Young and Y. Choi
London Math. Soc. Lect. Notes, 347 (2008) 248--342. arXiv:math/0604045

[Sk] А. Скопенков, Алгебраическая топология с элементарной точки зрения,
МЦНМО, Москва, в печати,
arXiv:math/0808.1395

[Th81] C.~Thomassen, Kuratowski's theorem, J.~Graph. Theory, 5 (1981) 225--242.

[Wh33] H. Whitney, Planar graphs, Fund. Math. 21 (1933) 73--84.

\normalsize